\documentclass{amsart}

\usepackage{amssymb}
\usepackage{amsthm}
\usepackage{amsmath}
\usepackage{amsbsy}
\usepackage[all]{xy}
\usepackage{bm}
\usepackage{hyperref}
\usepackage{tikz}
\usepackage{array}
\usepackage{colortbl,xcolor}
\newtheorem*{theorem}{Theorem}

\begin{document}

\title[]{Refining the Two-Dimensional \\Signed Small Ball Inequality} 

\author[]{Noah Kravitz}
\address[]{Grace Hopper College, Yale University, New Haven, CT 06510, USA}
\email{noah.kravitz@yale.edu}

\begin{abstract} 
The two-dimensional signed small ball inequality states that for all possible choices of signs,
$$ \left\| \sum_{|R| = 2^{-n}}{ \varepsilon_R h_R} \right\|_{L^{\infty}} \gtrsim n,$$
where the summation runs over all dyadic rectangles in the unit square and $h_R$ denotes the associated Haar
function.  This inequality first appeared in the work of Talagrand, and alternative proofs are due to Temlyakov and Bilyk \& Feldheim (who showed that the supremum equals $n+1$ in all cases).  We prove a stronger result: for all integers $0\leq k \leq n+1$, all possible choices of signs, and all dyadic rectangles $Q$ with $|Q| \geq 2^{-n-1}$,
$$ \left| \left\{ x \in Q: \sum_{|R| = 2^{-n}}{ \varepsilon_R h_R} = n + 1 - 2k\right\} \right| = \frac{|Q|}{2^{n+1}}\binom{n+1}{k}.$$
\end{abstract}
\maketitle

\vspace{-0.5cm}

\section{Introduction and Main Result}
\subsection{Introduction.} Let $\mathcal{D}$ denote the set of all dyadic intervals (i.e., intervals of the form $[m2^{-k}, (m+1)2^{-k})$ for some $m,k \in \mathbb{N}$) on $[0,1)$.  On such an interval $I$, we consider the associated Haar function
$$ h_{I}(x) = -\textbf{1}_{\text{left}}(x) + \textbf{1}_{\text{right}}(x)$$
that assumes the value $-1$ on the left half of the interval and $+1$ on the right half of the interval.  Dyadic rectangles in $[0,1)^d$ are defined naturally via Cartesian
products, and if $I = I_1 \times I_2 \times \dots \times I_d$, then $$h_{I}(x_1, x_2, \dots x_d) = h_{I_1}(x_1) h_{I_2}(x_2) \dots h_{I_d}(x_d)$$ extends the
definition of the associated Haar function. The \textit{signed small ball conjecture} states that for all choices of sign $\varepsilon_R \in \left\{-1,1\right\}$,
$$ \left\| \sum_{|R| = 2^{-n}}{ \varepsilon_R h_R} \right\|_{L^{\infty}([0,1)^d)} \gtrsim n^{\frac{d}{2}},$$
where the sum runs over all elements of $\mathcal{D}$ of fixed volume. This problem for $d=2$ first appeared in the 1994 work of Talagrand \cite{tala}
in connection with properties of the Brownian sheet.  Temlyakov \cite{tem} gave a surprising proof by adapting a classical argument of Sidon \cite{sidon} for lacunary Fourier series to this novel setting.  It is now understood that this problem has connections to deep and
long-standing problems in geometric discrepancy theory \cite{bil1, bil2}.  A third proof for $d=2$ appeared
recently in the work of Bilyk \& Feldheim \cite{bilyk}, whose argument gives a slightly sharper result: regardless of the choice of signs, there are always $2^{n+1}$ cubes of area $2^{-2n -2}$ on which the function assumes the value $n+1$.  They also established the
first fully rigorous connection to geometric discrepancy theory by showing a connection to binary nets.  The problem in higher dimensions
is wide open. The best current results are due to Bilyk-Lacey \cite{bil3}, Bilyk-Lacey-Vagharshakyan \cite{bil4, bil5} and 
 Bilyk-Lacey-Parissis-Vagharshakyan \cite{bil6},  and a special case was analyzed by Karslidis \cite{kar}.

\subsection{Main Result.}The purpose of this short paper is to generalize the refinement of Bilyk \& Feldheim and to show that, in a certain sense, the 
two-dimensional signed small ball inequality is surprisingly rigid. 
\begin{theorem} For any choice of signs $\varepsilon_R$, any integer $0\leq k \leq n+1$, and any dyadic rectangle $Q\in \mathcal{D}$ with $|Q| \geq 2^{-n-1}$,
$$ \left| \left\{ x \in Q: \sum_{|R| = 2^{-n}}{ \varepsilon_R h_R} = n + 1 - 2k\right\} \right| = \frac{|Q|}{2^{n+1}}\binom{n+1}{k}.$$
\end{theorem}
That is, the level sets on any dyadic rectangle of size at least $2^{-n-1}$ follow a binomial distribution, and the maximum value $n+1$ is assumed on each such rectangle (so that the two-dimensional signed small ball inequality is an immediate corollary).  In particular, setting $Q=[0,1)^2$ gives
$$ \left| \left\{ x \in [0,1)^2: \sum_{|R| = 2^{-n}}{ \varepsilon_R h_R} = n + 1 - 2k\right\} \right| = \frac{1}{2^{n+1}}\binom{n+1}{k},$$
and this statement for $k=0$ recovers the result of Bilyk \& Feldheim.  We believe that geometric interpretations of these level sets (especially in terms of discrepancy) could lead to interesting results in the future.
\\

Unfortunately, this line of reasoning, like earlier arguments, does not seem to generalize to higher dimensions. 
(In fact, it is fairly easy to construct explicit examples for all $d\geq 3$ where different choices of signs lead to different level sets.)  A pessimistic interpretation is that even though it has been assumed that the signed small ball inequality is understood for $d=2$, our result indicates that the $d=2$ case is so rigid that calling it an inequality
is a bit of a stretch: for an actual inequality to occur, we require $d\geq 3$, and this case remains far from understood.  In general, one could define two choices of signs to be in the same equivalence class if the arising distributions of level sets are the same. A natural topic of further inquiry is how quickly the number of equivalence classes grows with $n$ for $d\geq 3$.
\\

We remark that the distribution of level sets approximates a Gaussian distribution with standard deviation $\sim \sqrt{n}$. This is clearly of independent
interest and might have an interpretation in terms of properties of the Brownian sheet, where this inequality also arises.
It easily follows that for any $0 < p < \infty$, there exists a constant $c_p$ such that
$$ \left\| \sum_{|R| = 2^{-n}}{ \varepsilon_R h_R} \right\|_{L^{p}([0,1)^2)} \sim c_p\sqrt{n} \qquad \mbox{as}~n \rightarrow \infty.$$
$L^2-$orthogonality yields $c_2=1$. The fact that $L^{\infty}$ behaves substantially differently from $L^p$ for $1 \leq p < \infty$ further
underscores that this problem is quite subtle.

\section{Proof of the Main Result}
We first fix definitions.  For a fixed $n\in \mathbb{N}$, each $R\in \mathcal{D}$ of area $2^{-n}$ must have length $2^{-k}$ and height $2^{k-n}$ for some integer $0 \leq k \leq n$.  It is clear that for each such $k$, the $2^n$ dyadic rectangles of dimension $2^{-k} \times 2^{k-n}$ partition the unit square.  We call this set of dyadic rectangles, specified by $k$, the \textit{layer} $\mathcal{L}_k$.\\

The main idea of the proof is that we can flip the sign of each rectangle with negative sign without changing the distribution of level sets on $Q$.  We first show that for low layers, flipping a sign does not affect the distribution of level sets with respect to Haar functions in higher layers.  If all the Haar functions at lower layers have $\varepsilon_R = 1$, then flipping the sign does not affect the distribution of level sets at lower layers, either.  The same goes for high layers.  This allows us (proceeding from $\mathcal{L}_0$ and $\mathcal{L}_n$ towards the middle) to flip negative signs without ever changing the level set distribution.  Another way to think of this construction is as follows: `glue' the value of the function onto each point, then cut the unit square into smaller squares each with side length $2^{-n-1}$ and rearrange these small squares to reach the configuration where every $\varepsilon_R=1$.

\begin{proof}
Fix some $n\geq 0$ and some assignment of the $\varepsilon_R$'s.  Let $Q$ have length $2^{-a}$ and height $2^{-b}$.  Consider a dyadic rectangle $R_1\in \mathcal{L}_k$ with $k \leq n-b$.  This rectangle cannot overlap with any other rectangle in $\mathcal{L}_k$, but it overlaps with rectangles in layers $\mathcal{L}_{\ell}$ with $\ell < k$ and layers $\mathcal{L}_m$ with $m > k$.  Suppose $R_2\in \mathcal{L}_m$ with $m>k$, and consider the (interesting only if nonempty) intersection
$$R_3 = R_1 \cap R_2 \cap Q.$$
 $R_3$ is also a dyadic rectangle.  Whenever two dyadic intervals intersect, one is fully contained in the other.  Thus, by considering the heights of $R_1$ and $R_2$, we see that $R_3$ is completely contained in either the top or bottom half of $R_2$.  Note that within a single half of $R_2 = R_{2,1} \times R_{2,2}$, the function $h_{R_2}(x,y) = h_{R_{2,1}}(x) h_{R_{2,2}}(y)$ is independent of $y$, so $h_{R_2}(x,y)$ is independent of $y$ on $R_3$. Thus, the functions
$$h_{R_1} + \sum_{R \in \mathcal{L}_{m} \atop m > k}{ \varepsilon_R h_R} \quad \mbox{and} \quad  -h_{R_1} + \sum_{R \in \mathcal{L}_{m} \atop m > k}{ \varepsilon_R h_R}$$
have the same level set distribution on $R_1 \cap Q$ and hence also on $Q$ (since $k \leq n-b$, i.e., $Q$ is at least as tall as $R_1$). This shows that flipping the sign of $R_1$ does not affect the distribution of values with respect to higher layers. In particular, we may immediately set the signs of all $R \in \mathcal{L}_0$ to $\varepsilon_R=1$ (as long as $b\leq n$; see below for the case $b=n+1$).
Let us now assume that $\varepsilon_{R}=1$ for all $R\in \mathcal{L}_{\ell}$ for all $\ell < k\leq n-b$. It is easy to see that for any $\ell < k$, the function
$$\sum_{R\in \mathcal{L}_{\ell}}{ h_R}$$
is identical on the top and bottom halves of $R_1 \cap Q$.
So we see that the functions
$$h_{R_1} + \sum_{|R|=2^{-n} \atop R \neq R_1}{\varepsilon_R h_R} \quad \text{and} \quad -h_{R_1}+\sum_{|R|=2^{-n} \atop R \neq R_1}{\varepsilon_R h_R}$$

have the same level set distribution on $R_1 \cap Q$ and hence also on $Q$ as long as $\varepsilon_R = 1$ for all $R \in \mathcal{L}_{\ell}$ with $\ell < k$.
We can thus set $\varepsilon_{R}=1$ for all $R\in \mathcal{L}_k$ proceeding from $\mathcal{L}_0$ to $\mathcal{L}_{n-b}$.
By the same reasoning, we can also flip the signs of the rectangles in layers $\mathcal{L}_s$ with $a\leq s \leq n$ starting from $\mathcal{L}_n$ and working our way down to $\mathcal{L}_{a}$.  Recall that $|Q|=2^{-a-b} \geq 2^{-n-1} \rightarrow a+b \leq n+1 \rightarrow a \leq (n-b)+1$, so these two processes ``meet'' in the middle.  (In the extremal case $b=n+1$ mentioned above, all the work happens ``from the top''.)  As such, we can flip all signs to $\varepsilon_R=1$ without changing the level set distribution on $Q$.
An explicit computation in the case where $\varepsilon_R=1$ for all $R$ shows that
$$ \left| \left\{ x \in Q: \sum_{|R| = 2^{-n}}{ \varepsilon_R h_R} = n + 1 - 2k\right\} \right| = \frac{|Q|}{2^{n+1}}\binom{n+1}{k}.$$
\end{proof}

\textbf{Remarks.}  \begin{enumerate}
\item The requirement $|Q| \geq 2^{-n-1}$ is tight.  For smaller dyadic rectangles $Q$, the two processes do not necessarily meet in the middle, and for any $n$ it is easy to find such rectangles that do not exhibit the desired distribution.
\item Considering any tiling of the unit square by dyadic rectangles of area $2^{-n-1}$ yields a new geometric interpretation of the result of Bilyk \& Feldheim.
\end{enumerate}
\vspace{.3cm}

\textbf{Acknowledgment.} The author is grateful to Stefan Steinerberger for proposing this problem and engaging in helpful discussions throughout the writing process.

\end{document}